\documentclass[12pt,a4paper]{article}
\usepackage[latin1]{inputenc}
\usepackage{amsmath,amssymb}

\newcommand{\Legendre}[4]{\left(\frac{#1}{#2}\right)_{#3}^{#4}}
\newcommand{\mltsum}[2]{\sum\limits_{\stackrel{\scriptstyle #1}{#2}}}

\author{S.J.Patterson}
\title{The Fourier coefficients of metaplectic \\theta series on GL(2) over rational function fields
\thanks{The results of this papers were presented at a conference held in the Perrotis Agricultural School 
in Thessaloniki from 14$^\mathrm{th}$ to 18$^\mathrm{th}$ July, 2014 to celebrate Jeff Hoffstein's 60th birthday.   I would like 
to thank the organisers for the opportunity to speak at it. }}
\date{}
\begin{document}
\maketitle
\begin{center}To Jeff Hoffstein, to mark his 60$^\mathrm{th}$ birthday. \end{center}
\medskip
\begin{abstract} The idea of the metaplectic theta function was introduced by Tomio Kubota in the 1960s.
These theta functions are constructed as residues of Eisenstein series and are only known completely in 
the case of double covers and, up to the ambiguity of a constant, for triple covers.   In 1992 Jeff Hoffstein
published formul\ae\ by which these can be computed in certain cases over a rational function field.   The 
author gave an alternative approach in 2007.  Both of these methods give the coefficients in a closed form.  
The rational function field is unusual in that it has a large automorphism group.  In this paper we show that 
this group has an operation on the coefficients.   This operation is not visible from the explicit formul\ae.  
\end{abstract}

\section{Introduction}
The purpose of this paper is to take up the investigation begin in 
the paper \cite{P3}.    The main result of this paper is Theorem 1
which describes the action of the group of automorphisms of the
rational function field on the coefficients of metaplectic theta series.

We fix an integer $n$ and an odd prime power $q$, $q\equiv 1 \pmod n$.  The 
coefficients of the theta functions are quantities denoted by 
$\rho_0(r,\varepsilon,i)$ where $i$ is an integer  $0\le i<n$ and $r$ is a 
non--zero element of $\mathbb F_q[x]$.  The quantities describe the 
asymptotic distribution of Gauss sums of order $n$ over  $\mathbb F_q[x]$.
A crucial property of $\rho_0(r,\varepsilon,i)$ is that it depends only on 
$r$ modulo $n^\mathrm{th}$ powers.   Although the $\rho_0(r,\varepsilon,i)$ were
constructed only for polynomial $r$ and they can be extended to arbitrary 
rational functions $r$ by this property.

The automorphism group of $\mathbb F_q(x)$ over $\mathbb F_q$ is $\mathrm{PGL}_2(\mathbb F_q)$.
Let $\begin{pmatrix}a&b \cr c &d\end{pmatrix}$ be a non-degenerate matrix over 
$\mathbb F_q$.   We shall show that $\rho_0(r_1,\varepsilon,i)=\rho_0(r_2,\varepsilon,i)$ 
where $r_2(x)=r_1((ax+b)/(cx+d))((ad-bc)/(cx+d)^2)^{1-i}$.   It follows that there are 
a large number of relationships between the  $\rho_0(r,\varepsilon,i)$.   In particular
this very much simplifies the task of computing tables of this function for given
$n$ and $q$.

Section 2 is dedicated to summarizing and completing the results of \cite{P3} and in 
Section 3 we give the proof of Theorem 1.   It should be noted that the fact that 
such a statement should hold follows from the general discussion of the formalism
of the $\rho_0$ given in \cite{P1} and \cite{P2}.

In Section 4 we shall discuss the consequences of Theorem 1 when there are
at most three irreducible factors of $r$ and these are of degree $1$.  
Although these are very special they exhibit a number of interesting features
which are suggestive.   In particular these coefficients are, as was already 
noted in \cite{P3}, related to Selberg sums and some generalizations of these.   
Selberg sums were evaluated by G.W. Anderson, \cite{GWA} , R.J. Evans, \cite{RE} 
and P.B. van Wamelen, \cite{vW}, and the theory of the metaplectic group throws a 
new light on the general class of these sums.

\section{Summary of notations and previous results} We shall here both establish 
the notations we shall need and recall thoseresults from \cite{P3} which we shall 
make use of here.   The notations will,  as the result of experience differ a little 
from those in the previous paper. 

We let $q$ and $n$ be as above and let $k=\mathbb F_q(x)$ and $R=\mathbb F_q[x]$.   
We define the map $\chi: \mathbb F_q^\times \rightarrow \mu_n(\mathbb F_q); 
x \mapsto x^{(q-1)/n}$ where for any field $k$ $\mu_n(k)$ is the set of 
$n^\mathrm{th}$ roots of $1$ in $k$.   Let $\varepsilon$ be an embedding of 
$\mu_n(\mathbb F_q)$ in $\mathbb C^\times$.   Let 
$e_o:F_q \rightarrow \mathbb C^\times$ be a additive character.   It is convenient, and 
involves no loss of generality for our purposes, to assume that $e_o$ maps $a$ to 
$e^{2\pi \imath j/p}$ where the residue class $j \pmod p$ represents 
$\mathrm{Tr}_{\mathbb F_q/\mathbb F_p}(a)$.  We define an additive character $e$ on $k$ by 
$e(f)= e_o(\sum_v \mathrm{Res}_v(f  d x)) = e_o(-\mathrm{Res}_\infty(f \mathrm d x))$ 
where the sum over $v$ is over all finite places of $k$.   Note that in terms of the 
uniformizer  $x_\infty = x^{-1}$ at $\infty$ the latter expression is 
$e_o(\mathrm{Res}_\infty(f x_\infty^{-2}  \mathrm d x_\infty)$.   We define the Gauss 
sums over $R$ to be 
$g(r,\varepsilon,c)= \sum_{\xi \pmod c} \varepsilon(\Legendre{\xi}{c}{n}{})e_o(r\xi/c).$
Here $r$ and $c$ are non-zero elements of $R$.   For a character $\omega$ of 
$\mathbb F_q^\times$ we define the finite field Gauss sum $\tau(\omega) =
\sum_{j \in \mathbb F_q} \omega(j)e_o(j)$.   

The Davenport-Hasse theorem implies that for $r$ coprime to $c$ one has
\[
g(r,\varepsilon,c)=\mu(c) \varepsilon(\Legendre{r}{c}{n}{-1}
\Legendre{c'}{c}{n}{})(-\tau(\varepsilon\chi))^{\deg(c)}
\]
where $\mu$ denotes the M\"obius function in $R$. The case where $c$ and $r$ are 
no longer assumed to be coprime can be reduced to this case.   We shall come back 
to this later.

The functions which concern us here are 
\[
\psi(r,\varepsilon,\eta,s) = (1-q^{n-ns})^{-1} 
\sum_{c \in R,c \sim \eta} g(r,\varepsilon,c) q^{-\deg(c)s}
\]
where $s \in \mathbb C, \mathrm{Re}(s)> 3/2$ and 
$\eta \in k_\infty^\times/k_\infty^{\times n}$
where $k_\infty$ denotes the completion of $k$ at the infinite place.  
The condition $c \sim \eta$ means that $\eta/c \in R^\times k_\infty^{\times n}$.
The function $\psi(r,\varepsilon,\eta,s)$ has at most one pole
modulo $\frac {2\pi \sqrt{-1}}{n \log q}\mathbb Z$ in 
$\mathrm{Re}(s)> 1$ located at $s = 1 + 1/n$.   If it exists it is simple
and we denote by $\rho(r,\varepsilon,\eta)$ the residue of 
$\psi(r,\varepsilon,\eta,s)$ at $s = 1 + 1/n$.

The general theory of Eisenstein series over function fields (see
\cite{KP}) also shows that there exists a polynomial 
$\Psi(r,\varepsilon,i,T)$ so that 
\[
\psi(r,\varepsilon,\pi_\infty^{-i},s)=q^{-is}(1-q^{n+1-ns})^{-1}
\Psi(r,\varepsilon,i,q^{-ns}).
\]
The function $\Psi(r,\varepsilon,i,T)$ depends only on $i \pmod n$.
This leads to 
\[
\rho(r,\varepsilon,\pi_\infty^{-i})=c_1 q^{-i(n+1)/n} \Psi(r,\varepsilon,i,q^{-n-1})
\]
where $c_1=1/(n \log q)$.   The notation in \cite{P3} is different; there the residue
was replaced by the value of $(1-q^{1+\frac 1n -s})\psi(r,\varepsilon,\pi_\infty^{-i},s)$
at $s=1+1/n$.   This is inessential for our purposes.   We shall write  
\[
\rho_0(r,\varepsilon,i)=\Psi(r,\varepsilon,i,q^{-n-1})
\]
which is a much more convenient function to use.

We have also for $i$ with $0\le i <n$
\[
(1-q^{n-ns})\psi(r,\varepsilon,\pi_\infty^{-i},s)=\frac{q-1}n
\mltsum{i'\ge i}{i'\equiv i \pmod n} C(r,\varepsilon,i')q^{-i's}
\]
where
\[
C(r,\varepsilon,i) =  
\mltsum{\deg(c)=i}{c\ \mathrm{monic}} g(r,\varepsilon,c).
\]
If we let
\[
C^*(r,\varepsilon,i) =  
\mltsum{\deg(c)=i}{c\ \mathrm{monic},\mathrm{gcd}(r,c)=1} g(r,\varepsilon,c)
\]
then 
\[
C(r,\varepsilon,i)=\sum_{r^*}  g(r,\varepsilon,r^*)
\varepsilon(\chi(-1))^{(i-1)\deg(r^*)}
C^*(r{r^*}^{(n-2)},\varepsilon,i-\deg(r^*))
\]
where $r^*$ runs through the set of integers (modulo units) all of whose 
prime factors divide $r$ and where $g(r,\varepsilon,r^*)\neq 0$.  If 
we assume, as we shall do in this paper, that no non-trivial 
$(n-1)^{\mathrm{st}}$ power divides $r$ then we can describe the set of 
$r^*$ as follows.   Let $\Sigma$ be the set of primes dividing $r$, each of which
is to be represented by the corresponding monic polynomial.   Denote, for 
$\pi \in S$ the exponent of $\pi$ dividing $r$ by $e(\pi)$.   Then the 
set of $r^*$ can be parametrized by the subsets $S \subset \Sigma$ via
$r^*(S) = \prod_{\pi \in S} \pi ^{e(\pi)+1}$ with the restriction
$\sum_{\pi \in S} (e(\pi)+1)\mathrm{deg}(\pi) \le i$.   One can reduce the 
$g(r,\varepsilon,r^*)$ to a product of Legendre symbols, Gauss sums
of the type $\tau(\varepsilon\chi^j)$ and powers of $q$.   The general formula 
is neither illuminating nor computationally helpful so that we shall not discuss 
it here.   If we are involved in calculating $C(r,\varepsilon,i)$ and 
$i \le e(\pi)\mathrm{deg}(\pi)$ for all $\pi | r$ then the only $r^*$ in the
sum is $1$.

If we use the Davenport-Hasse theorem to evaluate the Gauss sums in the 
definition of $C^*(r,\varepsilon,i)$ then we see that 
\[
C^*(r,\varepsilon,i)=(-1)^i \tau(\varepsilon\chi)^i\sum_c \mu(c)
\overline{\varepsilon(\left(\frac{r}{c} \right)_n)}
\varepsilon(\left(\frac{c'}{c} \right)_n)
\]
where the sum over $c$ is as before.
The polynomial $\Psi(r,\varepsilon,i,T)$ satisfies a functional equation 
which is given implicitly in \cite[top of p. 251]{P3}.   One can make this much 
more useful by a simple observation.   Let $\sigma = \mathrm{deg}(r)+1$
and let $i$ be so that $0\le i < n$.   Let $R= [(\sigma - i)/n]$   Let 
$i'$ be the least non-negative residue modulo $n$ of $\sigma - i$.   We 
note first that $[(\sigma - i')/n]=R$.   To see this we observe that 
we obtain an equivalent statement if we replace $\sigma$ by $\sigma_1$ 
where $\sigma \equiv \sigma_1 \pmod n$.   We can then assume that 
$i \le \sigma_1 < i+n$ so that $R=0$.   Then $i'=\sigma_1 - i$ and 
$\sigma_1 - i' =i$.   As $0\le i < n$ we see that $[(\sigma_1 - i')/n]=0$
as required.

Next let $R_1=[(\sigma - 2i)/n]$ and $R_2=[(\sigma - 2i')/n]$.  Assume
that $i \neq i'$.  It is clear that $R-1 \le R_1,R_2 \le R$.  Assume that 
$\sigma_1$ is as above; then $\sigma_1 - 2i'= 2i-\sigma_1$ and so 
$[(\sigma_1 - 2i)/n]+[(\sigma_1 - 2i')/n]=-1$.   In the general case we therefore
have $R_1+R_2 =2 R - 1$.   Thus if $i<i'$ we have $R_1=R$ and $R_2=R-1$.

With these observations we can now formulate the functional equation 
for $\Psi(r,\varepsilon,i,*)$ in two equivalent forms:
\[
\begin{array}{rcl}
\Psi(r,\varepsilon,i,T)& =& (Tp^n)^{R+1} \frac{1-p^{-1}}{1-p^{n-1}T}\Psi(r,\varepsilon,i,(p^{2n}T)^{-1}
)\\
&& +\varepsilon\chi(-1)^{i\sigma}\tau(\varepsilon\chi^{i-i'})(p^nT)^R\Psi(r,\varepsilon,i',(p^{2n}T)^{-1})
\end{array} 
\]
and 
\[
\begin{array}{rcl}
\Psi(r,\varepsilon,i',T)& =& (Tq^n)^R \frac{1-q^{-1}}{1-q^{n-1}T}\Psi(r,\varepsilon,i',(q^{2n}T)^{-1}
)\\
&& +\varepsilon\chi(-1)^{i'\sigma}\tau(\varepsilon\chi^{i'-i})(q^nT)^R\Psi(r,\varepsilon,i,(q^{2n}T)^{-1}).
\end{array}
\]
Here $\sigma$ and $R$ are as above and we assume here and henceforth that $i < i'$.

With these notations the ``Hecke relations at infinity'' can be formulated as follows:
\[
\rho_0(r,\varepsilon,i) = \varepsilon\chi(-1)^{i\sigma}\tau(\varepsilon\chi^{i-i'})q^{-1}\rho_0(r,\varepsilon,i')
\]
and 
\[
\rho_0(r,\varepsilon,i') = \varepsilon\chi(-1)^{i'\sigma}\tau(\varepsilon\chi^{i'-i})\rho_0(r,\varepsilon,i).
\]
If $i=i'$ then $\rho_0(r,\varepsilon,i) = 0$ which is the reason for excluding this case
in in the functional equations.

The corresponding ``Hecke relation'' for a finite prime $\pi$ takes the form for $r_o$ 
coprime to $\pi$:
\[
\begin{array}{rl}
\rho_0(r_o\pi^j,\varepsilon,i)=&\varepsilon\chi(-1)^{i(j+1)\deg(\pi)}q^{((i)_n-(i-(j+1)\mathrm{deg}(\pi))_n-((j+1)\mathrm{deg}(\pi))_n)(1+\frac 1n)}\\
&q^{((j+1)\mathrm{deg}(\pi))_n}q^{-\left[ \frac{(j+1)\mathrm{deg}(\pi)}n\right]} g(-r_o,\varepsilon^{j+1},\pi)\\
&\rho_0(r_o\pi^{n-2-j},\varepsilon,i-(j+1)\deg(\pi))
\end{array}
\]
for $0\le j < n-1$ and 
\[
\rho_0(r_o\pi^{n-1},\varepsilon,i)=0.
\]
In \cite{P3} the relation (7) refers not to $\rho(r,\varepsilon,i)$ as used here but 
to the corresponding residue of $\psi(r,\varepsilon,i)$ so that the power of 
$\mathrm N(\pi)$ has had to be modified correspondingly.   It is convenient to
avoid fractional powers of $q$.

We retain the assumption that $i<i'$.   Because of the functional equations for 
the $\Psi(r,\varepsilon,i,T)$ we can determine the polynomial from the coefficients 
of $\Psi(r,\varepsilon,i,T)$ up to that of $T^{R/2}$ together with those of 
$\Psi(r,\varepsilon,i',T)$ up to $T^{R/2 - 1}$ if $R$ is even; if $R$ is odd these 
are both to be replaced by $T^{(R-1)/2}$.   This means that we need only evaluate 
the corresponding $C(r,\varepsilon,j)$ We shall give more details of the cases 
in Section 4.

We should stress here that although the $\rho(r,\varepsilon,i)$ only depend on $r$ modulo
$n^{\mathrm{th}}$ powers the same is not true of $\Psi(r,\varepsilon,i,T)$.  

We saw above that the global Gauss sums for $k$ could be expressed 
explicitly in terms of the Gauss sums for the field of constants and certain other
quantities.   First of all the connection between the Legendre symbol and the 
resultant shows that 
\[
\Legendre{c'}{c}{n}{} = \left\{\begin{array}{lr}
\chi(D(c))&\mbox{      if    } \mathrm{deg}(c) \equiv 0, 1 \pmod 4 \\
\chi(-D(c))&\mbox{      if    } \mathrm{deg}(c) \equiv 2, 3 \pmod 4
\end{array}\right.
\] 
where $D(c)$ is the discriminant of the polynomial $c$ (assumed monic).

The evaluation of the quadratic Gauss sum is also valid in $k$; details are given
in, for example, \cite[XIII,\S 12]{BNT}.    This yields a relation, known as  
Pellet's formula, 
\[
\omega(D(c)) = \mu(c)(-1)^{\mathrm{deg}(c)},
\]
valid of monic polynomials $c$ where $\omega$ denotes the quadratic character
of $\mathbb F_q^\times$.  Several proofs are available - see, for example, \cite{RGS}.   
It is this formula that establishes the connection of the $C^*(r,\varepsilon,i)$ in the 
special case $r=x^{e_0}(x-1)^{e_1}$ with a Selberg sum.

\section{The main theorem}

In \cite[\S 3]{P3} we explained how an analogue of $\psi(r,\varepsilon,\eta,s)$ can be 
defined over a ring obtained from $R$ through the inversion of an arbitrary finite set 
of primes and the relationship of these new funcctions with the original ones.
The formul\ae\ given there need a little further explanation.   The variable ``$\eta$''
in $\psi_{S}(r,\varepsilon,\eta,s)$ and $\psi_{S\cup \{\pi\}}(r,\varepsilon,\eta,s)$ should 
strictly speaking be considered as elements of $k_S^\times$ and $k_S^\times\times k_\pi^\times$
respectively.   What was written as $\psi_{S\cup\{\pi\}}(r,\varepsilon,\eta,s)$ should be
$\sum_{\theta \in r_\pi^\times/{\pi^\times n}}\psi_{S\cup\{\pi\}}(r,\varepsilon,\eta \times \theta,s)$;
here $r_\pi$ denotes the ring of integers of $k_\pi$, the completion of $k$ at $\pi$.   
 
Now we can move to the following theorem:

\noindent\textbf{Theorem 1} 
\emph{Let $g = \left( \begin{matrix} a&b\\c&d \end{matrix} \right)\in \mathrm{GL}_2(\mathbb F_q)$
and set $\Delta =\det(g)$.  Then if $r \in k$ 
and $r_i^g$ is defined by $r_i^g(x)=r(\frac{ax+b}{cx+d})\left(\frac{\Delta}{(cx+d)^2}\right)^{1-i}$ then}
\[
\rho(r,\varepsilon,i) = \rho(r_i^g,\varepsilon,i)
\]
\noindent\emph{Proof:} The case where $c=0$ is elementary so that we can concentrate
on the case $c\neq 0$.   Let $\pi_1 = x+d/c$ and $\pi_2 = x-a/c$.   Let $R_1$ (resp.
$R_2$) be the ring obtained from $R$ by inverting $\pi_1$ (resp. $\pi_2$).   Consider the 
map $g:k\rightarrow k; f(x) \mapsto f^g(x)=f((ax+b)/(cx+d))$.  Set $\pi_\infty = x^{-1}$.   
Then $\pi_\infty^g=-c^2\pi_1/\Delta + \mathcal{O}(\pi_1^2)$ (in $k_{\pi_1}$) and 
$\pi_2^g= -(\Delta/c^2) \pi_\infty +\mathcal{O}(\pi_\infty^2)$ (in $k_\infty$).   This means 
that $g$ maps $R_2$ to $R_1$.   We also have $\mathrm d (x^g)=\Delta \mathrm d x/(cx+d)^2$.
We can conclude that, for $\theta \in \mathbb F_q^\times$,
\[
\psi_{\{\infty,a/c\}}(r,\varepsilon,\pi_\infty^{m_1}\times \theta\pi_2^{m_2},s)
\] is
equal to 
\[
\psi_{\{-d/c,\infty\}}(r^g(x) (\Delta/(cx+d)^2),\varepsilon,(\delta^{-1}\pi_1)^{m_1}\times \theta(\delta\pi_\infty)^{m_2},s)
\]
where $\delta = -\Delta/c^2$.   Using the transformation properties of $\psi_{\{\infty,a/c\}}$ 
(resp. $\psi_{\{-d/c,\infty\}}$) under units of $R_2$ (resp. $R_1$) we deduce that 
\[
\psi_{\{\infty,a/c\}}(r,\varepsilon,\pi_\infty^{m_1}\pi_2^{-m_2}\times \theta,s)
\varepsilon\chi(\theta)^{-m_2}\varepsilon\chi(-1)^{m_2+m_1m_2}
\]
is equal to 
\[
\begin{array}{l}
\psi_{\{-d/c,\infty\}}(r^g(x)(\Delta/(cx+d)^2),\varepsilon,\theta^{-1}\times
(\delta^{-1}\pi_1)^{-m_1}(\delta\pi_\infty)^{m_2},s)\\
\hspace{4cm}\times\varepsilon\chi(\theta)^{2m_1+m_2}\varepsilon\chi(-1)^{m_1+m_1m_2}\varepsilon\chi(\delta)^{(m_1+m_2)^2};
\end{array}
\]
consequently 
\[
\psi_{\{\infty,a/c\}}(r,\varepsilon,\pi_\infty^{m_1}\pi_2^{-m_2}\times \theta,s)
\]
is equal to 
\[
\begin{array}{l}
\psi_{\{-d/c,\infty\}}(r^g(x)(\Delta/(cx+d)^2),\varepsilon,\theta^{-1}\times
(\delta\pi_\infty)^{m_1+m_2},s)\\
\hspace{4cm}\times\varepsilon\chi(\theta)^{2m_1+2m_2}\varepsilon\chi(-1)^{m_1+m_2}\varepsilon\chi(\delta)^{m_1+m_2}.
\end{array}
\]
Again the behaviour of $\psi_{\{-d/c,\infty\}}$ in the first variable under multiplication
by a unit (in our case this will be $((cx+d)^2/\Delta)^{m_1+m_2}$ yields that
\[
\psi_{\{-d/c,\infty\}}(((cx+d)^2/\Delta)^{-m_1-m_2}r^g(x)(\Delta/(cx+d)^2),\varepsilon,\theta^{-1}\times
(\delta\pi_\infty)^{m_1+m_2},s)
\]
is equal to 
\[
\psi_{\{-d/c,\infty\}}(r^g(x)(\Delta/(cx+d)^2),\varepsilon,\theta^{-1}\times
(\delta\pi_\infty)^{m_1+m_2},s)\varepsilon\chi(\theta)^{2(m_1+m_2)}\varepsilon\chi(-1)^{m_1+m_2}
\varepsilon\chi(\delta)^{(m_1+m_2)^2}.
\]
Together these show that 
\[
\psi_{\{\infty,a/c\}}(r,\varepsilon,\pi_\infty^{m_1}\pi_2^{-m_2}\times \theta,s)
\]
is equal to 
\[
\psi_{\{-d/c,\infty\}}(((cx+d)^2/\Delta)^{-m_1-m_2-1}r^g(x),\varepsilon,\theta^{-1}\times
(\delta\pi_\infty)^{m_1+m_2},s).
\]
We take the residue at $s=1+1/n$ and sum over $\theta \in \mathbb F_q^\times$.   We 
now get 
\[
\rho_{\{\infty\}}(r,\varepsilon,\pi_\infty^{m_1}\pi_2^{-m_2})=
\rho_{\{\infty\}}(((cx+d)^2/\Delta)^{-m_1-m_2-1}r^g(x),\varepsilon,
(\delta \pi_\infty)^{m_1+m_2}).
\]

The behaviour under units shows that both sides depend on $\delta$ in the 
same way.  We therefore obtain the formula of the theorem if we take 
residues and write $i$ for $-m_1-m_2$.

\section{Consequences}

We noted in \cite{P3} the elementary statement
\[
\rho_0(r,\varepsilon,i)\in \tau(\varepsilon\chi^i) \times  
q^{-(n+1)[(1+\deg(r)-i)/n]}\mathbb Z[1^{1/n}]
\]
where the estimate for the denominator can probably be improved.   The real challenge 
is to find the factor lying in $\mathbb Z[1^{1/n}]$.   It follows from the results of 
Section 2 that $\rho_0(r,\varepsilon,i)=0$ if $R=[(1+\deg(r)-i)/n]<0$.   Let $i$ and 
$i'$ be as in \S2 with $i<i'$ and $0\le i,i' <1$.  Recall that we have 
$R=[(1+\deg(r)-i')/n]$.   The pair $\rho_0(r,\varepsilon,i)$ and $\rho_0(r,\varepsilon,i')$ 
are connected with one another by the Hecke relation at infinity.  In
this connection we recall that $\tau(\varepsilon\chi^a)^i$ is 
$(-1)^{i-1}\tau(\varepsilon\chi^{ai})$ times an (integral) Jacobi sum.

We shall in this section consider the cases where $r$ has at most three prime 
factors and that these are of degree $1$.   By means of a linear transformation 
we can assume then that $r$ has the form $x^{e_0}(x-1)^{e_1}(x-\lambda)^{e_\lambda}$
where $\lambda \in \mathbb F_q$ and $\lambda \neq 0,1$.   By the Hecke relations at 
$x$, $x-1$ and $x-\lambda$ we can take $e_j\le [n/2]-1$ for $j=0,1,\lambda$. 

If just one prime divides we can take $r$ to be  $x^{e_0}$, the degree of which is 
$e_0$.    It follows that that $\rho_0(x^{e_0},\varepsilon,i)=0$ if $i>e_0+1$.   By 
means of the Hecke relations at infinity we only need to determine these coefficients 
when $i < (e_0+1)/2$.   This means that $C(x^{e_0},\varepsilon,i)=C^*(x^{e_0},\varepsilon,i)$.
This sum is a degenerate Selberg sum which we could evaluate by means of the results
in \cite{vW}.

We can however proceed in a different way which is, in some respects, illuminating.   
We apply Theorem 1 with $\left(\begin{matrix} 0 & -1 \\ 1 & 0\end{matrix}\right)$ and 
we see that $\rho_0(x^{e_0},\varepsilon,i)=\rho_0(x^{-e_0+2i-2},\varepsilon,i)$   We
can now apply the Hecke relation at $x$ and we see that this latter function is 
a multiple of $\rho_0(x^{e_0-2i},\varepsilon,e_0+1-i)$; the multiplier is made up 
of 
a power of $q$, a power of $\varepsilon\chi(-1)$ and $\tau(\varepsilon\chi^{2i-e_0-1})$. 
We note that $i=0$ then $\rho_0(x^{e_0},\varepsilon,i) =1$.   If $i>0$ then $e_0-2i>-1$ 
and so $0\le e_0-2i <e_0$.   This means that we can compute
$\rho_0(x^{e_0},\varepsilon,i)$ recursively.   We deduce that
$\rho_0(x^{e_0},\varepsilon,i)$ is $\tau(\varepsilon\chi^i)$ times a Jacobi sum, that is, a
product of elementary Jacobi sums.   This is also what the theory of Selberg sums yields.

We now move on to the case where there are two prime factors, thus $r$ is of the form 
$x^{e_0}(x-1)^{e_1}$ where $1\le e_0,e_1 \le [n/2]-1$.   Now $R=[(e_0+e_1+1-i)/n]$ and 
we see that $\rho_0(x^{e_0}(x-1)^{e_1},\varepsilon,i) = 0$ if $i>e_0+e_1+1$.   Again
the Hecke relation at infinity means that we can restrict our attention to the case
$i \le [(e_0+e_1-1)/2]$.   It is clear that we can exchange the roles of $e_0$ and 
$e_1$ without changing the function.   We may assume also that $e_0\le e_1$.
As we have already studied the case where one of the exponents is $0$ we assume that 
$e_0 \ge 1$.

With these restrictions we have $R=0$ in the notation used above.   This means 
that, as we shall see below, $\rho_0(x^{e_0}(x-1)^{e_1},\varepsilon,i)$ is equal to 
$c_1C(x^{e_0}(x-1)^{e_1},\varepsilon,i)$ where now $i$ is to be understood not as 
a residue class but as the least non-negative element of that class.   There are 
two cases to be distinguished.   If $e_o+1>i$ then 
$C(x^{e_0}(x-1)^{e_1},\varepsilon,i)=C^*(x^{e_0}(x-1)^{e_1},\varepsilon,i)$
and this means that $\rho_0(x^{e_0}(x-1)^{e_1},\varepsilon,i)$ is equal to
$\tau(\varepsilon\chi^i)$ times a Selberg sum which itself is a Jacobi sum.

If $e_o+1 \le i$ then there are precisely two elements $r^*$ in the terminology of 
Section 2, $1$ and $x^{e_0+1}$.   We can now apply Theorem 1 again but now inverting 
$x-1$ then we find that $\rho_0(x^{e_0}(x-1)^{e_1},\varepsilon,i) = 
\rho_0((x-1)^{2i-2-e_0-e_1}x^{e_0},\varepsilon,i)$.   If we now make use of the Hecke 
relation at $(x-1)$ we obtain a multiple of 
$\rho_0((x-1)^{e_0+e_1-2i}x^{e_0},\varepsilon,e_0+e_1-i+1)$.   The case $i=0$ does not 
arise and so $e_0+e_1-2i +1 < e_0+e_1-i+1$ whence it follows that 
$\rho_0((x-1)^{e_0+e_1-2i}x^{e_0},\varepsilon,e_0+e_1-i+1)=0$; indeed the corresponding 
$\Psi$ function vanishes.   It therefore follows that if $e_o+1 \le i$
\[
\rho_0(x^{e_0}(x-1)^{e_1},\varepsilon,i) =0
\]
and in the other case, namely $e_o+1 > i$ 
\[
\rho_0(x^{e_0}(x-1)^{e_1},\varepsilon,i)
\]
is a power of $q$ times $c_1\tau(\varepsilon\chi^i)$ times a specific Jacobi sum.   We shall come 
back to the details in a later publication.   What is rather remarkable is that there is,
in these cases, an ``explicit formula'' for $\rho_0(x^{e_0}(x-1)^{e_1},\varepsilon,i)$ and
that this is in terms of Jacobi sums.   I have not found a method of demonstrating this
without the use of the theorem of Anderson, Evans and v. Wamelen.

This is not the rule.   If we next consider $r$ of the form 
$x^{e_0}(x-1)^{e_1}(x-\lambda)^{e_\lambda}$, now with $0< e_0,e_1,e_\lambda \le [n/2]-1$,
then we find a large number of relationships between various
$\rho_0(x^{e_0}(x-1)^{e_1}(x-\lambda)^{e_\lambda},\varepsilon,i)$.
The structure of this set of relations is that which one knows from the theory of the 
hypergeometric function (see \cite{W&W}).   If we have $\mathrm{Min}(e_0,e_1,e_\lambda)+1>i$ 
and $e_0+e_1+e_2 \le n$ then we see, as before, that 
$\rho_0(x^{e_0}(x-1)^{e_1}(x-\lambda)^{e_\lambda},\varepsilon,i)
=C^*(x^{e_0}(x-1)^{e_1}(x-\lambda)^{e_\lambda},\varepsilon,i)$.   The latter sum is 
an analogue of the hypergeometric function in the same sense that the standard 
Selberg sum is an analogue of the beta function.

In certain special cases, when $n=4$, or when $n=6$ and $e_0,e_1,e_\lambda$ are all $2$ 
conjectures of Eckhardt and Patterson and of Chinta, Friedburg and Hoffstein respectively
suggest that these sums take on a special form and this is what one finds.  It seems 
as if the corresponding statements can now be proved by a new method due to S. Friedberg
and D.Ginzburg \cite{FG}  but there is still work to be done to complete the proof.     
In other cases, as one would expect, the values are irregular.    These evaluations go 
beyond the range considered in \cite{JH1} and make it clear that with increasing complexity 
of the $r$ the arithmetical nature of the $\rho_0(r,\varepsilon,i)$ also becomes more complex.
Furthermore a complete evaluation does not seem to reasonable expectation and one will 
have to be satisfied with less specific questions.    One can, for example, make estimates
for the $\rho_0(r,\varepsilon,i)$ in different metrics.   One would suspect that one can 
do better than relatively elementary convexity bounds.

We now return to the expression of $\rho_0(r,\varepsilon,i)$ in terms of  $C(r,\varepsilon,i)$ 
or of $C^*(r,\varepsilon,i)$ for $R=0,1,2$ where, as before, $R=[(1+\deg(r)-i)/n]$.   Let 
$0\le i <n$ and $i'=(1+\deg(r)-i)_n$.    We assume, as before that $i < i'$   The condition 
$R \le 3$ means that we cover all cases with $\deg(r) < 3n$.   This is much more than we
need for the purposes of this paper.    

We recall that by construction 
\[
\Psi(r,\varepsilon,i,T)=
\frac{1-q^{n+1}T}{1-q^{n}T}\frac{q-1}n
\mltsum{i'\ge i}{i'\equiv i \pmod n} C(r,\varepsilon,i')T^{(i'-i)/n}.
\]
We shall assume that $0\le i <n$. 

If $R=0$ then $\Psi(r,\varepsilon,i,T)=C(r,\varepsilon,i)$ and there
is nothing left to be said.

If $R=1$ then a direct application of the definition yields 
\[
\Psi(r,\varepsilon,i,T)=C(r,\varepsilon,i)+(C(r,\varepsilon,i+n)
-(q-1)q^nC(r,\varepsilon,i))T.
\] 

If $R=2$ then
\[
\begin{array}{r}
\Psi(r,\varepsilon,i,T)=C(r,\varepsilon,i)+(C(r,\varepsilon,i+n)
-(q-1)q^nC(r,\varepsilon,i))T \\+(C(r,\varepsilon,i+2n)
-(q-1)q^nC(r,\varepsilon,i+n)-(q-1)q^{2n}C(r,\varepsilon,i))T^2
\end{array}
\]

It follows that in the three cases we have 
\[
\Psi(r,\varepsilon,i,q^{-n-1})= C(r,\varepsilon,i),
\]
\[
\Psi(r,\varepsilon,i,q^{-n-1})=q^{-1}C(r,\varepsilon,i)+q^{-n-1}C(r,\varepsilon,i+n)
\] 
and 
\[
\begin{array}{r}
\Psi(r,\varepsilon,i,q^{-n-1})=C(r,\varepsilon,i)q^{-2}+
C(r,\varepsilon,i+n)q^{-n-2}\\+C(r,\varepsilon,i+2n)q^{-n-2}.
\end{array}
\]
respectively.   These are perfectly usable expressions but are not as practical as they 
might be as the number of summands in $C(r,\varepsilon,i+n)$ and $C(r,\varepsilon,i+2n)$ 
is large.  It is advantageous to exploit the functional equation when $R=1$ and $R=2$; 
the method is a simple version of the ``approximate function equation''.   

To carry this out we write $C_j = C(r,\varepsilon,i+jn)$ and $C'_j = C(r,\varepsilon,i'+jn)$.
We shall use $X$ with $X=q^nT$ rather than $T$. Let $F(X)=\Psi(r,\varepsilon,i,q^{-n}X)$ and
$G(X)=\Psi(r,\varepsilon,i',q^{-n}X)$.    Then $F(X) = \sum_{0\le j \le R} D_j X^j$ and 
$G(X) = \sum_{0\le j \le R} D_j' X^j$ where 
$D_j=C_jq^{-nj} -(q-1)\sum_{0\le \ell < j}C_\ell q^{-n\ell}$ and
$D_j'=C_j'q^{-nj} -(q-1)\sum_{0\le \ell < j}C_{\ell}' q^{-n\ell}$.   Let 
$\eta = \varepsilon\chi(-1)^{i\sigma}\tau(\varepsilon\chi^{i-i'})$ and 
$\eta' = \varepsilon\chi(-1)^{i\sigma}\tau(\varepsilon\chi^{i'-i})$.   Note 
that $\eta \eta' = q^{-1}$.   Then the functional equation takes on the 
two equivalent forms
\[
F(X) = x^{R+1} \frac{1-q^{-1}}{1-q^{-1}X}F(\frac 1X)+
\eta \frac{1-X}{1-q^{-1}X}G(\frac 1X)
\]
and 
\[
G(X) = x^{R} \frac{1-q^{-1}}{1-q^{-1}X}G(\frac 1X)+
\eta' \frac{1-X}{1-q^{-1}X}F(\frac 1X)
\]
We can rewrite the first as 
\[
G(X) = \eta^{-1}\left( X^RF(\frac 1X) +(1-q^{-1}) \frac{F(X)-X^RF(\frac 1X)}{1-X} \right)
\]
and the second as 
\[
F(X) = \eta'^{-1}\left( X^RG(\frac 1X) +(1-q^{-1}) \frac{XG(X)-X^RG(\frac 1X)}{1-X} \right).
\]
These can be written as linear expression for the $D_*$ in terms of the $D_*'$ or conversely.
Explicitly one has
\[
D_k= \eta'^{-1}(D_{R-k}' -(1-q^{-1})D_k'+(1-q^{-1})\sum_{j<\min (k,R-k)} D_j' -D_{R-j-1}')
\]
and 
\[
D_k'= \eta^{-1}(D_{R-k} -(1-q^{-1})\sum_{j<\min (k+1,R-k)} D_j -D_{R-j-1}).
\]
These are not precisely what we need for an ``approximate functional equation''.   What is 
need to do is to express the $D_*$ and $D_*'$ in terms of $D_0,\ldots D_{[R/2]}$ and 
$D_0',\ldots D_{[R/2]}'$.   This is easy to do for any given value of $R$.  Thus we find
if $R=0$ 
\[
D_0=\eta D_0',
\]
if $R=1$ 
\[
D_1= \eta'^{-1}D_0' - (q-1)D_0
\]
and if $R=2$
\[
D_2= \eta'^{-1}D_0' - (q-1)D_0.
\]
It follows that in the three cases we have $F(X)=D_0$, $F(X)=D_0(1- (q-1)X)+\eta'^{-1}D_0'X$
and $F(X)=D_0(1- (q-1)X^2)+D_1X+\eta'^{-1}D_0'X^2$ respectively.   We obtain the values of
$F(1/q)$ in terms of $D_0$, of $D_0$ and $D_0'$ and of $D_0$, $D_1$ and $D_0'$ respectively.

It should be noted that whereas on the one hand the theory of Selberg sums allows us 
to obtain closed expressions for the coefficients of metaplectic sums on the other hand the 
theory of metaplectic forms, and, in particular Theorem 1, leads to a number of new relations 
between Selberg sums which do not seem to be accessible by elementary methods or those of 
Anderson.   This will be the subject of a future paper.

\section{Outlook - Curves of higher genus}
The case discussed above is very straightforward as the structure of the
rational curve is explicit.   For other curves there is, in general, no 
``natural'' ring of integers, especially if we demand that it be a principal 
ideal ring.   It seems therefore a considerable challenge to gather numerical 
evidence in such a case.   It is unclear as to whether the nature of the 
$\rho_0(r,\varepsilon,i)$ are typical of what happens in the case of curves of
genus $\ge 1$.   There seems to be no reason why not but nevertheless one is 
able to exploit so much in the rational case that one is cautious about 
making any too large extrappolations.   

One case that is probably deserving of study is that of elliptic curves. 
To gain some idea of what is needed we consider the function field of an 
elliptic curve over a field of characteristic $\neq 2,3$.   We can then 
represent it in Weierstrass form and as usual we make the point at infinity 
the identity element of the group.   Hasse's estimate shows that there is 
at least one further rational point on the curve and one knows that the group 
of points over a finite field is either a (non-trivial) cyclic group or a
product of two cyclic groups; see \cite[Prop. 7.1.9]{HC}.   Let $P_1$ (resp. $P_1$ amd $P_2$) be 
generators of the group of rational points in these two cases.   Then the 
ring $R$ of functions integral outside $\{\infty,P_1\}$ (resp. $\{\infty,P_1,P_2\}$)
is a principal ideal domain.   The methods of the theory of algebraic
curves allow us to represent the ring $R$ as a quotient of a polynomial ring.   
This comes down to working with a model of the curve in 3- or 4-dimensional
projective space.   The determination of the elements of this ring and especially 
of the prime elements is possible but it is no longer as easy as in the case of 
the rational function field.   Also the description of the relation $\sim$ 
in $k_S^{\times}$ is considerably more intricate than before.   It seems at the 
moment that the computational effort needed would be considerable and that
it would only be justified it there were good reason to expect that the 
behaviour of the $\psi(r,\varepsilon,\eta,s)$ and $\rho_0(r,\varepsilon,\eta)$
would show features that do not appear in the case of rational function fields.
Whether this is so is an open question at present.

\noindent
Mathematisches Institut\\
Bunsenstr. 3--5\\
37073 G\"ottingen\\
Germany

\noindent
e-mail:\texttt{spatter@gwdg.de}

\end{document}